\documentclass[12pt]{article}
\usepackage{mathrsfs}
\usepackage{epic,eepic,epsf,epsfig}
\usepackage{amsfonts,srcltx,mathrsfs}
\usepackage{multirow}
\usepackage{amsmath}
\usepackage{amssymb}
\usepackage{amsbsy}
\usepackage{graphicx}
\usepackage{amsfonts}
\usepackage{color}
\newtheorem{lem}{Lemma}[section]%
\newtheorem{theorem}[lem]{Theorem}%
\newtheorem{cor}[lem]{Corollary}%
\newtheorem{exam}[lem]{Example}%
\newtheorem{prop}[lem]{Proposition}%

\def\a{\alpha}

 \def\O{\Omega} \def\G{\Gamma}

\def\di{\bigm|}  
\def\nd{\mathrel{\bigm|\kern-.7em/}}

\def\f{\noindent}

\def\PSL{\hbox{\rm PSL}}\def\PSU{\hbox{\rm PSU}}
 \def\PGL{\hbox{\rm PGL}} \def\GL{\hbox{\rm GL}} \def\Mult{\hbox{\rm Mult}}
\def\PSp{\hbox{\rm PSp}}\def\P\GammaL{\hbox{\rm P\Gamma L}} 
 \def\ASL{\hbox{\rm ASL}}

\def\Aut{\hbox{\rm Aut}}

\def\Syl{\hbox{\rm Syl}}

\def\Cos{\hbox{\rm Cos}}

\newcommand{\qed}{\mbox{\raisebox{0.7ex}{\fbox{}}} \vspace{4truemm}}
\def\mz{{\mathbb Z}}

\begin{document}
\title{Pentavalent symmetric graphs admitting vertex-transitive non-abelian simple groups}

\author{Jia-Li Du, Yan-Quan Feng\footnotemark, Jin-Xin Zhou\\
{\small\em Department of Mathematics, Beijing
Jiaotong University, Beijing 100044, China}}

\footnotetext[1]{Corresponding author. E-mails:
JiaLiDu@bjtu.edu.cn, yqfeng@bjtu.edu.cn, jxzhou@bjtu.edu.cn}

\date{}
 \maketitle

\begin{abstract}

A graph $\G$ is said to be {\em symmetric} if its automorphism group $\Aut(\G)$ is transitive on the arc set of $\G$. Let $G$ be  a finite non-abelian simple group and let $\G$ be a connected pentavalent symmetric graph such that $G\leq \Aut(\G)$. In this paper, we show that if $G$ is transitive on the vertex set of $\G$, then either $G\unlhd \Aut(\G)$ or $\Aut(\G)$ contains a non-abelian simple normal subgroup $T$ such that $G\leq T$ and $(G,T)$ is one of $58$ possible pairs of non-abelian simple groups. In particular, if $G$ is arc-transitive, then $(G,T)$ is one of $17$ possible pairs, and if $G$ is regular on the vertex set of $\G$, then $(G,T)$ is one of $13$ possible pairs, which improves the result on pentavalent symmetric Cayley graph given by Fang, Ma and Wang in 2011.

\bigskip
\f {\bf Keywords:} Symmetric graph, Cayley graph, coset graph, simple group.

\medskip
\f {\bf 2010 Mathematics Subject Classification:} 05C25, 20B25.

\end{abstract}

\section{Introduction}

Let $G$ be a permutation group on a set $\O$ and let $\a\in \O$. Denote by $G_{\a}$ the stabilizer of $\a$ in $G$, that is, the subgroup of $G$ fixing the point $\a$. We say that $G$ is {\em regular} on $\O$ if for any two points there is a unique element of $G$ mapping one to the other.
Denote by $\mz_n$, $D_n$, $A_n$ and $S_n$ the cyclic group of order $n$, the dihedral group of order $2n$, the alternating group and the symmetric group of degree $n$, respectively. For a subgroup $H$ of a group $G$, denote by $C_G(H)$ the centralizer of $H$ in $G$ and by $N_G(H)$ the normalizer of $H$ in $G$.

Throughout this paper, all groups and graphs are finite, and all graphs are simple and undirected. For a graph $\G$, we denote its vertex set and automorphism group by $V(\G)$ and $\Aut(\G)$, respectively. A graph $\G$ is said to be {\it $G$-vertex-transitive} for $G\leq \Aut(\Gamma)$ if $G$ acts transitively on $V(\G)$, {\it $G$-regular} if $G$ acts regularly on $V(\G)$, and {\it $G$-symmetric} if $G$ acts transitively on the arc set of $\G$ (an arc is an ordered pair of adjacent vertices). In particular, $\G$ is {\it vertex-transitive} or {\it symmetric} if it is $\Aut(\G)$-vertex-transitive or $\Aut(\G)$-symmetric, respectively. The graph $\G$ is a {\em Cayley graph} on $G$ if $G$ is regular on the vertex set of
$\G$, and the Cayley graph is \emph{normal} if $G$ is a normal subgroup of $\Aut(\G)$.

Let $G$ be a non-abelian simple group and let $\G$ be a pentavalent symmetric $G$-vertex-transitive graph. In this paper, we show that either $G\unlhd \Aut(\G)$ or $\Aut(\G)$ contains a non-abelian simple normal subgroup $T$ such that $G<T$ and $(G,T)$ is one of $58$ possible pairs of non-abelian simple groups. The motivation of this investigation comes from the two extreme cases, that is, $G$-symmetric and $G$-regular.

There are two steps to study a symmetric graph $\G$ --- the first step is to investigate the normal quotient graph for a normal subgroup of an arc-transitive group of automorphisms (see Section~\ref{s2} for the definition of quotient graph), and the second step is to reconstruct
the original graph $\G$ from the normal quotient by using covering
techniques. This is usually done by taking a maximal normal subgroup such that the normal quotient graph has the same valency as the graph $\G$, and in this case, the normal quotient graph is called a {\em basic graph} of $\G$. The situation seems to be somewhat more promising with $2$-arc-transitive graphs (a $2$-arc is a directed path of length $2$), and the strategy for the structural analysis of these graphs, based on taking normal quotients, was first laid out by Praeger (see \cite{P1,P3}). The strategy works for locally primitive graphs, that is, vertex-transitive graphs with vertex stabilizers acting primitively on the corresponding neighbors sets (see \cite{P4,P5}). For more results, refer to \cite{FengLZ,LiLZ} for example. For a non-abelian simple group $G$, the $G$-symmetric graphs are important basic graphs and have received wide attention. For example, Fang and Praeger~\cite{FangP,FP2} classified $G$-symmetric graphs admitting $G$ as Suzuki simple groups or Ree simple groups acting transitively on the set of $2$-arcs of the graphs. For cubic $G$-symmetric graph, it was proved by Li~\cite{CHLi} that either $G$ is a normal in $\Aut(\G)$, or $(G,\Aut(\G))=(A_7,A_8)$, $(A_7,S_8)$, $(A_7,2.A_8)$, $(A_{15},A_{16})$ or $(\GL(4,2),\rm AGL(4,2))$. Fang {\em et al}~\cite{FangLW} proved that none of the above five pairs can happen, that is, $G$ is always normal in $\Aut(\G)$. In this paper, we show that if $\G$ be a connected pentavalent $G$-symmetric graph, then either $G\unlhd \Aut(\G)$ or $\Aut(\G)$ contains a non-abelian simple normal subgroup $T$ such that $G<T$ and $(G,T)$ is one of $17$ possible pairs of non-abelian simple groups.

Investigation of Cayley graphs on a non-abelian simple group is currently a hot topic in algebraic graph theory. One of the most remarkable achievements is the complete classification of connected trivalent symmetric non-normal Cayley graphs on finite non-abelian simple groups. This work was began in 1996 by Li~\cite{CHLi}, and he proved that a connected trivalent symmetric Cayley graph $\G$ on a non-abelian simple group $G$ is either normal or $G=A_5$, $A_7$, $\PSL(2,11)$, $M_{11}$, $A_{11}$, $A_{15}$, $M_{23}$, $A_{23}$ or $A_{47}$. In 2005, Xu {\em et al}~\cite{XFWX2005} proved that either $\G$ is normal or $G= A_{47}$, and two years later, Xu {\em et al}~~\cite{XFWX} further showed that if $G=A_{47}$ and $\G$ is not normal, then $\G$ must be $5$-arc-transitive and up to isomorphism there are exactly two such graphs.

Let $\G$ be a connected symmetric Cayley graph $\G$ on a non-abelian simple group $G$. Fang, Praeger and Wang~\cite{FangP2} developed a general method to investigate the automorphism group $\Aut(\G)$, and using this, the symmetry of Cayley graphs of valency $4$ or $5$ has been investigated in~\cite{Fang4,Fang5}. Now let $\G$ be non-normal and of valency $5$. Fang~{\em et al}~\cite{Fang5} proved that if $\Aut(\G)$ is quasiprimitive then $(G,\mathrm{soc}(\rm Aut(\G)))=(A_{n-1}, A_n)$, where either $n=60\cdot k$ with $k\di  2^{15}\cdot3$ and $k \neq3, 4, 6, 8$ or $n = 10\cdot m$ with $m\di8$, and if $\Aut(\G)$ is not quasiprimitive then there is a maximal intransitive normal subgroup $K$ of $\Aut(\G)$  such that the socle of $\Aut(\G)/K $, denoted by $\bar{L}$, is a simple group containing $\bar{G}=GK/K\cong G$ properly, where $(\bar{G},\bar{L})=(\O^-(8,2),\rm PSp(8,2))$, $(A_{14},A_{16})$, $(A_{n-1},A_n)$ with $n\geq 6$ such that $n$ is a divisor of $2^{17}\cdot 3^2\cdot 5$ or $\bar{G}=\bar{L}$ is isomorphic to an irreducible subgroup of $\PSL(d,2)$ for $4\leq d \leq 17$ and the $2$-part $|G|_2>2^d$. In this paper, we prove that $\Aut(\G)$ contains a non-abelian simple normal subgroup $T$ such that $G<T$ and $(G,T)$ is one of $13$ possible pairs of non-abelian simple groups. The following is the main result of this paper.

\begin{theorem}\label{theo=main}
Let $G$ be a non-abelian simple group and $\G$ a connected pentavalent symmetric $G$-vertex-transitive graph. Then either $G\unlhd \Aut(\G)$, or $\Aut(\G)$ contains a non-abelian simple normal subgroup $T$ such that $G<T$ and $(G,T)=(\Omega^-_8(2),\PSp(8,2))$, $(A_{14},A_{16})$, $(\PSL(2,8),A_{9})$ or $(A_{n-1},A_n)$ with $n\geq 6$ and $n\di 2^9\cdot3^2\cdot5$.
\end{theorem}

\begin{cor}\label{cor=arc}
Let $G$ be a non-abelian simple group and $\G$ a connected pentavalent $G$-symmetric graph. Then either $G\unlhd \Aut(\G)$, or $\Aut(\G)$ contains a non-abelian simple normal subgroup $T$ such that $G<T$ and $(G,T)=(A_{n-1},A_n)$ with $n=2\cdot 3$, $2^2\cdot 3$, $2^4$, $2^3\cdot 3$, $2^5$, $2^2\cdot 3^2$, $2^4\cdot 3$, $2^3\cdot 3^2$, $2^5\cdot 3$, $2^4\cdot 3^2$,
$2^6\cdot 3$,  $2^5\cdot 3^2$, $2^7\cdot 3$, $2^6\cdot 3^2$, $2^7\cdot 3^2$,
$2^8\cdot 3^2$ or $2^9\cdot 3^2$.
Moreover, if $T=A_6$, then $\Aut(\G)=S_6$ and if $T=A_{2^5}$, $A_{2^7 \cdot 3}$, $A_{2^7 \cdot 3^2}$, $A_{2^8 \cdot 3^2}$ or $A_{2^9\cdot 3^2}$, then $\Aut(\G)=T$.
\end{cor}

\begin{cor}\label{cor=regular}
Let $G$ be a non-abelian simple group and $\G$ a connected pentavalent symmetric $G$-regular graph. Then either $G\unlhd \Aut(\G)$, or $\Aut(\G)$ contains a non-abelian simple normal subgroup $T$ such that $G<T$ and $(G,T)=(A_{n-1},A_n)$ with $n=2\cdot 5$, $2^2\cdot 5$, $2^3\cdot5$, $2\cdot3\cdot 5$, $2^4\cdot5$, $2^3\cdot 3\cdot5$, $2^4\cdot 3^2\cdot5$, $2^6\cdot 3\cdot5$, $2^5\cdot 3^2\cdot5$, $2^7\cdot 3\cdot5$,
$2^6\cdot 3^2\cdot5$,  $2^7\cdot 3^2\cdot5$ or $2^9\cdot 3^2\cdot5$.
\end{cor}

For connected  symmetric cubic Cayley graphs on  non-abelian simple groups, similar to Corollary~\ref{cor=regular} there are six possible pairs $(G,T)=(A_{47},A_{48})$, $(\PSL(2,11),M_{11})$, $(M_{11},M_{12})$, $(A_{11}, A_{12})$, $(M_{23},M_{24})$ or $(A_{23},A_{24})$ (see \cite[Theorem 7.1.3]{CHLi}), and Xu {\em et al}~\cite{XFWX2005,XFWX} proved that only the pair $(G,T)=(A_{47},A_{48})$ can happen and there are exactly two connected non-normal symmetric cubic Cayley graphs on $A_{47}$ with automorphism groups $A_{48}$. Based on the method given in~\cite{XFWX2005,XFWX}, we may show that $(G,T)=(A_{9},A_{10})$ or $(A_{19},A_{20})$ cannot happen in Corollary~\ref{cor=regular}. For $(G,T)=(A_{39},A_{40})$, with MAGMA~\cite{magma} we find a connected non-normal pentavalent Cayley graph on $A_{39}$ (see Example~\ref{NonNormalExample}).

\section{Preliminaries\label{s2}}

In this section, we describe some preliminary results which will be used later. First we describe vertex stabilizers of connected pentavalent symmetric graphs.

\begin{prop}\label{prop=stabilizer} {\rm \cite[Theorem 1.1]{Guo}}
Let $\Gamma$ be a connected pentavalent $G$-symmetric graph with $v\in V(\Gamma)$. Then $G_v\cong \mathbb{Z}_5$, $D_{5}$, $D_{10}$, $F_{20}$, $F_{20}\times \mathbb{Z}_2$,
$ F_{20}\times \mathbb{Z}_4$, $A_5$, $S_5$, $A_4\times A_5$, $S_4\times S_5$, $(A_4\times A_5)\rtimes \mathbb{Z}_2$, $ {\rm ASL}(2,4)$, ${\rm AGL}(2,4)$,
${\rm A\Sigma L}(2,4)$, ${\rm A\Gamma L}(2,4)$ or $\mathbb{Z}^6_2\rtimes {\rm \Gamma L}(2,4)$, where $F_{20}$ is the Frobenius group of order $20$,
 $A_4\rtimes \mathbb{Z}_2=S_4$ and $A_5\rtimes \mathbb{Z}_2=S_5$.
In particular, $|G_v|=5$, $2\cdot 5$, $2^2\cdot 5$, $2^2\cdot 5$, $2^3\cdot 5$, $2^4\cdot 5$, $2^2\cdot 3 \cdot 5$, $2^3\cdot 3 \cdot 5$,
$2^4\cdot 3^2\cdot 5$, $2^6\cdot 3^2\cdot 5$, $2^5\cdot 3^2\cdot 5$, $2^6\cdot3 \cdot 5$, $2^6\cdot 3^2\cdot 5$, $2^7\cdot 3\cdot 5$, $2^7\cdot 3^2\cdot 5$
or $2^9\cdot 3^2\cdot 5$, respectively.
\end{prop}

Let $\G$ be a graph and $N \leq \Aut(\G)$. The \emph{quotient graph} $\G_N$ of $\G$ relative to $N$ is defined as the graph with vertices the orbits of $N$ on $V(\G)$ and with two orbits adjacent if there is an edge in $\G$ between these two orbits.

\begin{prop}\label{prop=atlesst3orbits}{\rm \cite[Theorem 9]{Lorimer}}
Let $\Gamma$ be a connected $G$-symmetric graph of prime valency, and let $N\unlhd G$ have at least three orbits on $V(\Gamma)$. Then $N$ is the kernel of $G$ on $V(\G_N)$, and semiregular on $V(\G)$. Furthermore, $\G_N$
is $G/N$-symmetric with $G/N \leq \Aut(\G_N)$.
\end{prop}

For a group $G$ and a prime $p$, denote by $O_p(G)$ the largest normal $p$-subgroup of $G$ and by $\Phi(G)$ the Frattini subgroup of $G$, that is,
the intersection of all maximal subgroups of $G$. By \cite[Lemma 2.7]{Wang}, we have the following proposition.

\begin{prop}\label{largest normal subgroup}
For a group $G$ and a prime $p$, let $H=O_p(G)$ and $V=H/\Phi(H)$. Then $G$ has a natural action on $V$, induced from the action of $G$ on $H$ by conjugation.
If $C_G(H)\leq H$, then $H$ is the kernel of this action of $G$ on $V$, that is, $G/H\leq \GL(V)$.
\end{prop}

The following result follows from the classification of three-factor simple groups.

\begin{prop}\label{prop=235simplegroup}{\rm \cite[Theorem \MakeUppercase{\romannumeral1} ]{Huppert2}}
Let $G$ be a non-abelian simple $\{2,3,5\}$-group.
Then $G= A_5$, $A_6$ or $\rm PSU(4,2)$.
\end{prop}

Given a group $G$, its factorization $G=HD$ is said to be \emph{maximal} if both $H$ and $D$ are maximal subgroups of $G$.

\begin{prop}\label{prop=simplegroupinedx2a3b5c}{\rm \cite[Lemma 3.3]{Fang5}}
Let $T$ be a simple group and $G$ a non-abelian simple subgroup such that $|T:G|=n>1$ and $n=2^a\cdot 3^b\cdot 5^c$ for some $0\leq a\leq 9$, $0\leq b\leq 2$ and $0\leq c \leq1$. If $T=HD$ is a maximal factorization with $G\leq H$, then $T$, $G$ and $|T:G|$ are listed in Table~\ref{table=2}.
\end{prop}

\begin{table}[ht]

\begin{center}

\begin{tabular}{|c|c|c||c|c|c|}

\hline
$T$        & $G$           & $|T:G|$                                   & $T$              & $G$             & $|T:G|$ \\
\hline
$M_{11}$   & $\PSL(2,11)$  & $2^2\cdot3$                              &$M_{12}$          & $\PSL(2,11)$    & $2^4\cdot3^2$   \\
\hline
$M_{24}$   & $M_{23}$      & $2^3\cdot 3$                             & $M_{12}$         & $M_{11}$        & $2^2\cdot 3$  \\
\hline
$\rm P\Omega^{+}_8(2)$  & $\PSp(6,2)$   & $2^3\cdot 3\cdot 5$       & $\PSp(6,2)$      & $A_8$      & $2^3\cdot3^2$ \\
\hline
$\rm P\Omega^{+}_8(2)$  & $A_9$      & $2^6\cdot 3\cdot 5$             & $\PSp(6,2)$      & $A_7$      & $2^6\cdot3^2$ \\
\hline
$\PSp(8,2)$   & $\Omega^{-}_8(2)$  & $2^4\cdot 3\cdot 5$               & $\PSU(4,2)$    & $A_6$     & $2^3\cdot3^2$ \\
\hline
$\PSU(3,3)$       & $\PSL(2,7)$   & $2^{2}\cdot3^2$                   & $^2F_{4}(2)'$   & $\PSL(2,25)$  & $2^8\cdot3^2$\\
\hline
$\PSp(4,4)$       & $\PSL(2,16)$   & $2^4\cdot 3\cdot 5 $           & $\PSL(3,4)$     & $\PSL(2,7)$  & $2^3\cdot3\cdot 5$\\
\hline
$A_{16}$           & $A_{14}$      & $2^4\cdot 3\cdot 5$            &$A_{10}$           & $A_{8}$      & $2\cdot 3^2\cdot 5$  \\
\hline
$A_{10}$           & $A_{7}$      & $2^4\cdot 3^2\cdot 5$           & $A_{9}$           & $A_{7}$      & $2^3\cdot 3^2$  \\
\hline
$A_{9}$           & $\PSL(2,8)$      & $2^3\cdot 3^2\cdot 5$         &$A_{8}$           & $\PSL(3,2)$      & $2^3\cdot 3\cdot 5$  \\
\hline
$A_{7}$           & $\PSL(2,7)$      & $3\cdot 5$                   &$\rm PSp(6,2)$          & $\PSL(2,8)$  & $2^6\cdot3^2\cdot 5$ \\
\hline
$A_n$          & $A_{n-1}$  & $n=2^a\cdot3^b\cdot 5^c$   &&&\\
\hline
\end{tabular}

\end{center}
\vskip -0.5cm
\caption{{Non-abelian simple group pairs of index $2^a\cdot 3^b\cdot 5^c$}}\label{table=2}
\end{table}

Let $G$ and $E$ be two groups. We call an extension $E$ of $G$ by $N$ a {\em central extension} of $G$ if $E$ has a central subgroup $N$ such that $E/N\cong G$, and if further $E$ is perfect, that is, the derived group $E'$ equal to $E$, we call $E$ a {\em covering group} of $G$. A covering group $E$ of $G$ is called a {\em double cover} if $|E|=2|G|$. Schur~\cite{Schur} proved that for every non-abelian simple group $G$ there is a unique maximal covering group $M$ such that every covering group of $G$ is a factor group of $M$ (see \cite[Kapitel V, \S23]{Huppert}). This group $M$ is called the {\em full covering group} of $G$, and the center of $M$ is the {\em Schur multiplier} of $G$, denoted by $\Mult(G)$.

\begin{prop}\label{prop=covering group} For $n\geq 5$, the alternating group $A_n$ has a unique double cover $2.A_n$, and for $n\geq 7$,   all subgroups of index $n$ of $2.A_n$ are conjugate and isomorphic to $2.A_{n-1}$.
\end{prop}

\f {\bf Proof:} By Kleidman and Liebeck~\cite[Theorem 5.1.4]{Kleidman}, $\Mult(A_n)\cong \mz_2$ for $n\geq 5$ with $n\not=6,7$,  and $\Mult(A_n)\cong \mz_6$ for $n=6$ or $7$. This implies that $A_n$ has a unique double cover for $n\geq 5$, denoted by $2.A_n$. Since $A_n$ has no proper subgroup of index less than $n$, all subgroups of index $n$ of $2.A_n$ contain the center of $2.A_n$. Let $n\geq 7$. By~\cite[2.7.2]{Wilson}, $2.A_n$ contains a subgroup $2.A_{n-1}$ of index $n$, and since all subgroups of index $n$ of $A_n$ are conjugate, all subgroups of index $n$ of $2.A_n$ are conjugate and isomorphic to $2.A_{n-1}$.
\hfill\qed

Let $G$ be a group. For $H\leq G$, let $D$ be a union of some double cosets of $H$ in $G$ such that $D^{-1}=D$. The \emph{coset graph} $ \Gamma=\Cos(G,H,D)$ on $G$ with respect to $H$ and $D$ is defined to have vertex set $V(\Gamma)=[G:H]$, the set of right cosets of $H$ in $G$, and edge set $E(\Gamma)=\{\{Hg,Hdg\}\ |\ g\in G,d\in D\}$. The graph $\Gamma$ has valency $|D|/|H|$ and it is connected if and only if $G=\langle D,H\rangle$, that is, $D$ and $H$ generate $G$. The action of $G$ on $V(\G)$ by right multiplication induces a transitive group of automorphisms, and this group is symmetric if and only if $D$ is a single double coset. Moreover, this action is faithful if and only if $H_G=1$, where $H_G$ is the largest normal subgroup of $G$ contained in $H$.

Conversely, let $\Gamma$ be a $G$-vertex-transitive graph with $G\leq \Aut(\Gamma)$. By~\cite{Sabidussi}, the graph $\Gamma$ is isomorphic to a coset graph $\Cos(G,H,D)$, where $H=G_{v}$ is the vertex stabilizer of $v\in V(\Gamma)$ in $G$ and $D$ consists of all elements of $G$ mapping $v$ to one of its neighbors.
It is easy to show that $H_G=1$ and $D$ is a union of some double cosets of $H$ in $G$ satisfying $D^{-1}=D$. Assume that $\G$ is $G$-symmetric and $g \in G$ interchanges $v$ and one of its neighbors. Then $g^2 \in H$ and $D=HgH$. Furthermore, $g$ can be chosen as a $2$-element in $G$, and the valency of $\Gamma$ is $|D|/|H|=|H:H\cap H^g|$. For more details regarding coset graph, referee to  \cite{FangP,Lorimer,Miller,Sabidussi}.

\begin{prop}\label{prop=cosetgraph}
Let $\Gamma$ be a connected $G$-symmetric graph of valency $k$, and let $\{u,v\}$ be an edge of $\Gamma$. Then $ \Gamma $ is isomorphic to a coset graph $\Cos(G,G_{v},G_{v}gG_{v})$,
where $g$ is a $2$-element in $G$ such that $G_{uv}^{g}=G_{uv}$, $g^2\in G_{v}$,  $\langle G_v,g\rangle=G$, and $k=|G_v:G_v \cap G_{v}^{g}|$.
\end{prop}

In Proposition \ref{prop=cosetgraph}, the $2$-element $g$ such that $g^2\in G_{v}$,  $\langle G_v,g\rangle=G$ and $k=|G_v:G_v \cap G_{v}^{g}|$, is called {\em feasible} to $G$ and $G_v$. Feasible $g$ can be computed by MAGMA~\cite{magma} when the order $|G|$ is not too large, and this is often used in Section~\ref{theo=main}.

\section{Proof of Theorem~\ref{theo=main}}

In this section, we always assume that $G$ is a non-abelian simple group. Let us begin by proving a series
of lemmas which will be used in the proof of Theorem~\ref{theo=main}.

\begin{lem}\label{lem=A_5}
Let $\G$ a connected pentavalent symmetric $G$-vertex-transitive graph with $G=A_5$. Then either $G\unlhd \Aut(\G)$, or $\G=K_6$ and $\Aut(\G)=S_6$.
\end{lem}

\f {\bf Proof:} Set $A=\Aut(\G)$. Assume $G\ntrianglelefteq  A$. Since $|A_5|=2^2 \cdot 3\cdot 5$
and $\G$ is a connected pentavalent symmetric $G$-vertex-transitive graph, we have $
|V(\G)|=2^2 \cdot3 \cdot5$, $2 \cdot3 \cdot5$, $2^2 \cdot5$, $2 \cdot5 $, $2^2 \cdot3$ or $2 \cdot3$.
If $|V(\G)|=2^2 \cdot 3\cdot 5$, by \cite[Theorem 4.1]{Guo12p}, $\G=\mathcal{G}_{60}$ and $A=A_5 \rtimes D_5$,
implying that $G\unlhd A$, contrary to the assumption $G\ntrianglelefteq  A$. By \cite[Theorems~4.1 and 4.2]{hua},
there is no connected pentavalent symmetric graph of order $2\cdot 3\cdot 5$ or $2^2 \cdot 5$.
If $|V(\G)|=2\cdot 5$, then $\G=K_{5,5}$ with $A=(S_5 \times S_5)\rtimes \mz_2$ where
$S_5\times S_5$ fixes the bipartite sets of $\G$. Since $G=A_5\leq S_5\times S_5$,
$\G$ is not $G$-vertex-transitive, contrary to the hypothesis that $\G$ is $G$-vertex-transitive.
If $|V(\G)|=2^2 \cdot 3$, by \cite[Theorem 4.1]{hua}, $\G$ is the Icosahedron graph $I_{12}$,
or $K_{6,6}-6K_2$, the complete bipartite graph of order $12$ minus a one factor.
If $\G=K_{6,6}-6K_2$ then $A=S_6\times \mz_2$ and $S_6$ fixes the bipartite sets of $\G$.
Thus, $G=A_5\leq A_6\leq S_6$ and $\G$ is not $G$-vertex-transitive, a contradiction.
If $\G=I_{12}$ then $A=A_5 \times \mz_2$, implying that $G\unlhd A$, a contradiction.
Thus, $|V(\G)|=6$. It follows that $\G=K_6$ and $A= S_6$.
\hfill\qed

\begin{lem}\label{lem=GH}
Let $\G$ be a connected pentavalent symmetric $G$-vertex-transitive graph with $v\in V(\G)$. Let $H\leq \Aut(\G)$ and $GH\leq \Aut(\G)$.
Then $|(HG)_v|/|G_v|=|H|/|H\cap G|$ is a divisor of $2^9\cdot 3^2\cdot5$.
\end{lem}

\f {\bf Proof:} Since $\G$ is $G$-vertex-transitive, the Frattini argument implies $GH=G(GH)_v$ and
$|H||G|/|H \cap G|=|HG|=|G(GH)_v|=|G||(GH)_v|/|G\cap (GH)_v|=|G||(GH)_v|/|G_v|$,
that is, $|H|/|H\cap G|=|(GH)_v|/|G_v|$. Since $GH\leq \Aut(\G)$,  Proposition~\ref{prop=stabilizer} implies
$|(GH)_v|\di 2^9\cdot 3^2\cdot5$. It follows $|H|/|H\cap G|\di 2^9\cdot 3^2\cdot 5$.
\hfill\qed

The {\em radical} of a finite group is the largest solvable normal subgroup of the group.

\begin{lem}\label{lem=RG}
Let $\G$ be a connected pentavalent symmetric $G$-vertex-transitive graph and $R$ the radical of $\Aut(\G)$. Let $G\neq A_5$ and $5\di |R|$. Then $RG=R\times G$.
\end{lem}

\f {\bf Proof:} Let $v \in V(\G)$ and set $B=RG$. Since $G\cap R\unlhd G$, we have  $G\cap R=1$ because the simplicity of $G$ implies $G\not\leq R$. It follows $|B|=|R||G|$, and by Lemma~\ref{lem=GH}, $|R|=|B_v|/|G_v|\di 2^9\cdot 3^2\cdot 5$. Since $5 \di |R|$, we have $5\di |B_v|$, that is, $\G$ is $B$-symmetric.

Since $R$ is solvable, $R$ has a Hall $\{2,3\}$-subgroup, say $H$. Set $\Omega=\{H^r|r\in R\}$. By~\cite{Hall}, all Hall $\{2,3\}$-subgroups of $R$ are conjugate and so the conjugate action of $B$ is transitive on $\Omega$. Since $|R|\di 2^9\cdot 3^2\cdot 5$, we have $|\Omega|=|R:N_R(H)|=1$
or $5$, and since $G\neq A_5$, the conjugate action of $G$ on $\Omega$ is trivial. This means that $G$ normalizes $H$, implying $GH\leq B$. Set $Y=GH$ and $\Delta=\{Yb|b \in B\}$. Then $|\Delta|=|B:Y|=5$. The kernel of the right multiplication action of $B$ on $\Delta$ is $Y_B$, the largest normal subgroup of $B$ contained in $Y$. It follows $B/Y_B \lesssim S_5$ and so $G\leq Y_B$, as $G\neq A_5$.

Since $|G||H|/|G\cap H|=|GH|=|Y|=|G||Y_v|$, $Y_v$ is a $\{2,3\}$-group as $H$ is a $\{2,3\}$-group, and hence $(Y_B)_v$ is a $\{2,3\}$-group. Recall that $\G$ is $B$-symmetric and has valency $5$.
If $(Y_B)_v\not=1$, then $B_v$ is primitive on the neighborhood $\Gamma(v)$ of $v$ in $\G$, and then $Y_B\unlhd B$ implies $5\di |(Y_B)_v|$, a contradiction. It follows $(Y_B)_v=1$ and $Y_B=G(Y_B)_v=G$. In particular, $G\unlhd B$ and so $B=R\times G$.
\hfill\qed

\begin{lem}\label{lem=S_9}
Let $\G$ a connected pentavalent $G$-symmetric graph with $G\leq S_9$ and $G\neq A_5$. Then $G\unlhd \Aut(\G)$.
\end{lem}

\f {\bf Proof:} Set $A=\Aut(\G)$. By the simplicity of $G$, $G\leq S_9$ implies $G\leq A_9$. By the Atlas~\cite[p. 37]{Atlas}, $G=A_9$ or $G\leq A_8, A_7, A_6,A_5,\PSL(2,8)$.
Similarly, if $G\leq A_8$ then $G=A_8$ or $G\leq A_7, A_6, \PSL(3,2),A_5$;
if $G\leq A_7$ then $G=A_7$ or $G\leq A_6, A_5, \PSL(2,7)$; if $G\leq A_6$
then $G=A_6$ or $A_5$. Since $\G$ is $G$-symmetric with $G\neq A_5$, we have $G=A_9$, $A_8$, $A_7$ or $A_6$.

By Proposition~\ref{prop=cosetgraph}, $\G=\Cos(G,G_{v},G_{v}gG_{v})$ for some feasible $g$, that is, a $2$-element $g$ in $G$ satisfying $g^2\in G_{v}$,  $\langle G_v,g\rangle=G$ and $5=|G_v:G_v \cap G_{v}^{g}|$. By Proposition~\ref{prop=stabilizer}, $g$ has order at most $8$.

For $G= A_9$, by the Atlas \cite[p. 37]{Atlas} and Proposition~\ref{prop=stabilizer},
we have $G_v\cong \mz_5$, $D_5$, $D_{10}$, $F_{20}$, $F_{20}\times \mz_2$, $A_5$, $S_5$,
$A_4\times A_5$ or $(A_4\times A_5)\rtimes \mz_2$. By MAGMA, $G$ has two conjugacy classes of subgroups isomorphic to $ D_5, D_{10}, F_{20}, F_{20}\times \mz_2, A_5$ or $S_5$ and one conjugacy class of subgroups isomorphic to $\mz_5, A_4\times A_5$ or $(A_4\times A_5)\rtimes \mz_2$.

Let $G_v\cong A_5$ or $S_5$. By taking a given $G_v$ in each conjugacy class,
computation with MAGMA shows that there is no feasible $g$.

Let $G_v\cong D_5, D_{10}, F_{20}$ or $F_{20}\times \mz_2$. Then $G_v$ has two conjugacy classes,
and computation shows that there is no feasible $g$ in one of the two classes. Take a given $G_v$
in the other class. By MAGMA, if $G_v\cong D_5$ then there are $80$ feasible $g$, and for each of them,
the corresponding coset graph $\G$ is a connected pentavalent $G$-symmetric graph with $A=A_9\times \mz_2^2$.
Similarly, if $G_v\cong D_{10}$, $F_{20}$ or $F_{20}\times \mz_2$ then  in each case,
there are $80$ feasible $g$, and for each of them, the corresponding coset graph $\G$
is a connected pentavalent $G$-symmetric graph with $A=A_9\times \mz_2$, $A_9\rtimes \mz_4$ or $S_9$ respectively. It follows $G\unlhd A$.

Let $G_v\cong \mz_5, A_4 \times A_5$ or $(A_4 \times A_5)\rtimes \mz_2$. Then $G_v$ has one conjugacy
class and take a given $G_v$ in the class. By MAGMA, if $G_v\cong \mz_5$ then there are $120$ feasible $g$,
and for each of them, the corresponding coset graph is a connected pentavalent $G$-symmetric graph with $A= A_9\times \mz_4$. Similarly, if $G_v \cong A_4 \times A_5$ or $(A_4 \times A_5)\rtimes \mz_2$
then in each case, there are $480$ feasible $g$, and for each of them, the corresponding coset graph
is a connected pentavalent $G$-symmetric graph with $A= A_9\times \mz_2$ or $S_9$ respectively. It follows $G\unlhd A$.

For $G=A_8$, $A_7$ or $A_6$, similar arguments as the case $G=A_9$ give rise to $G\unlhd A$.
\hfill\qed

\begin{lem}\label{lem=GP} Let $\G$ be a connected pentavalent $G$-symmetric graph and $R$ the radical of $\Aut(\G)$.
Then $RG=R\times G$.
\end{lem}

\f {\bf Proof:} Set $B=RG$ and $A=\Aut(\G)$. Since $G\cap R\unlhd G$, we have  $G\cap R=1$ because the simplicity of $G$ implies $G\not\leq R$. It follows $|B|=|R||G|$.
Since $\G$ is $G$-symmetric, we have $5\di|G_v|$ for $v\in V(\G)$. By Lemma~\ref{lem=GH}, $|R|=|B_v|/|G_v|$, and by Proposition~\ref{prop=stabilizer}, $|R|\di 2^9\cdot3^2$.

If $G\leq A_9$, Lemmas \ref{lem=A_5} and \ref{lem=S_9} imply that either $G\unlhd A$, or $G=A_5$ and $\G=K_6$ (for this case $R=1$). Thus, $RG=R\times G$ for $G\leq A_9$. In what follows we may assume $G\nleq S_9$.

Let $\Syl_2(R)$ be the set of Sylow $2$-subgroups of $R$. Recall that $R\unlhd A$ and $|R|\di 2^9\cdot 3^2$.
By Sylow's Theorem, $|\Syl_2(R)|=1,3$ or $9$. Since $R\unlhd B$, $G$ has a natural action on $\Syl_2(R)$ by conjugation, and since $G\nleq S_9$, $G$ fixes every element in $\Syl_2(R)$. Thus, $G$ normalizes $P$ for any $P\in \Syl_2(R)$. Clearly, $|P|=2^r$ for some $0\leq r\leq 9$.

Now we claim $GP=G\times P$. Set $D=PG$. Then $P\unlhd D$ and $\G$ is $D$-symmetric. The simplicity of $G$ implies $P\cap G=1$. By Lemma~\ref{lem=GH}, $|P|=|P|/|P\cap G|=|D_v|/|G_v|=2^r$ and hence $|D_v|=2^r\cdot |G_v|$, which implies that if $G_v\cong \mz_5$ then by Proposition~\ref{prop=stabilizer}, $D_v\cong D_5$, $D_{10}$, $F_{20}$, $F_{20}\times \mz_2$ or $F_{20}\times \mz_4$. Thus, $r\leq 5$. This is always true for other $G_v$ by checking all other possible cases in  Proposition~\ref{prop=stabilizer}. In particular, $r=5$ if and only if $(G_v,D_v)\cong (A_5, {\rm A\Sigma L}(2,4))$ or
$(A_4\times A_5, \mathbb{Z}^6_2\rtimes  {\rm \Gamma L}(2,4))$.

Let $C=C_D(P)$. Then $C\unlhd D$ and $G\cap C=G$ or $1$. Suppose $G\cap C=1$. Since $P$ is $2$-group and $D/P\cong G$
is simple, we have $P=O_2(D)$, the largest normal $2$-subgroup of $D$. On the other hand, since $D=PG$ and $G\cap C=1$,
$C$ is $2$-group and hence $C\leq P$. Since $r\leq 5$, Proposition~\ref{largest normal subgroup} implies $G\cong D/P\leq \GL(5,2)\cong \PSL(5,2)$, and since $G\nleq S_9$, the Atlas \cite[p. 70]{Atlas} implies $G=\PSL(5,2)$ because $\PSL(4,2)\cong A_8$.
It follows that $r=5$ and $(G_v,D_v)\cong (A_5, {\rm A\Sigma L}(2,4))$ or $(A_4\times A_5, \mathbb{Z}^6_2\rtimes  {\rm \Gamma L}(2,4))$. However, by MAGMA, $G$ has no subgroup isomorphic to $A_4 \times A_5$. Thus, $G_v\cong A_5$, and by Proposition~\ref{prop=cosetgraph}, $\G=\Cos(G,G_v,G_vgG_v)$ for a $2$-element $g$ such that $g^2\in G_v$, $\langle G_v,g\rangle=G$ and $|G_v:G_v\cap G_v^g|=5$. Again by MAGMA, $G$ has four conjugacy classes of subgroups isomorphic to $A_5$ and by taking a $G_v$ in each class, computation shows that there is no feasible $g$ for $G_v$. It follows $G\cap C=G$, that is, $G\leq C$. Since $G\cap P=1$, we have $D=P\times G$, as claimed.

Set $\Omega=\{Db~|~b\in B\}$. Then $|\Omega|=1,3$ or $9$. By considering the right multiplication action of $B$ on $\Omega$,
we have $B/D_B \leq S_9$, where $D_B$ is the largest normal subgroup of $B$ contained in $D$. Since $G\nleq S_9$,
we have $G\leq D_B$. Thus, $G$ is characteristic in $D_B$ because $D=P\times G$ implies that $G$ is characteristic in $D$. It follows that $G\unlhd B$ as
$D_B\unlhd B$. Since $R\unlhd B$ and $R\cap G=1$, we have $B=RG=R \times G$.
\hfill\qed

\begin{lem}\label{lem=insolvable2}
Let $\G$ be a connected pentavalent $X$-symmetric $G$-vertex-transitive graph with $G\leq X$. Assume $X$ has trivial radical. Then either $G\unlhd X$, or $X$ has a non-abelian
simple normal subgroup $T$ such that $G<T$ and $(G,T)=(\Omega^-_8(2),\PSp(8,2))$, $(A_{14},A_{16})$, $(\PSL(2,8),A_{9})$ or $(A_{n-1},A_n)$ with $n\geq 6$ and $n\di 2^9\cdot3^2\cdot5$.
\end{lem}

\f {\bf Proof:} Let $A=\Aut(\G)$. Let $N$ be a minimal normal subgroup of $X$. Since $X$ has trivial radical,
$N=T^s$ for a positive integer $s$ and a non-abelian simple group $T$. Then $NG\leq X$, and by Lemma~\ref{lem=GH},
$|N|/|N\cap G|\di 2^9\cdot 3^2\cdot 5$.

If $G=A_5$, by Lemma \ref{lem=A_5}, either $G\unlhd X$ or $\G=K_6$ and $A=S_6$. For the latter, $(G,T)=(A_5,A_6)$. In what follows we may assume $G\neq A_5$ and $G\ntrianglelefteq X$. Thus, $G\ntrianglelefteq A$.

Suppose $N \cap G =1$. Then $|N|=|N|/|N\cap G|\di 2^9\cdot 3^2\cdot 5$, $|K|=|NG|=|N||G|$ and $K/N\cong G$. Write $K=NG$. For $v\in V(\G)$, we have $K=GK_v$ and $|K|=|GK_v|=|G||K_v|/|G_v|$. It follows $|K_v|=|G_v||N|$. Note that $|\PSU(4,2)|=2^6 \cdot3^4\cdot 5$. By Proposition~\ref{prop=235simplegroup}, $N=A_5$ or $A_6$.
Thus, $5\di |K_v|$, forcing that $\G$ is $K$-symmetric. Let $C=C_K(N)$.

Assume $C=1$. Then $G\leq K=K/C \leq \Aut(N)\leq S_6$, implying that $5\di |G|$ and $5^2 \nmid |K|$. On the other hand,
$5^2 \di |N||G|=|K|$, a contradiction.

Assume $ C\neq 1$. Then $C\cap N=Z(N)=1$, and $C=C/C\cap N\cong CN/N \unlhd K/N \cong G$. Since $|K|=|N||G|=|N||C|$,
we have $K=C\times N$. First let $C\cap G=1$. Then $|CG|=|C||G|$. Since $CG\leq K$, $|CG|$ is a divisor of $|K|=|N||C|$, implying $|G|\di |N|$. Since $N\cong A_5$ or $A_6$ and $G \neq A_5$, we have $C\cong G\cong N=A_6$ and hence
$K=C\times N\cong A_6\times A_6$. If $C_v=1$ then $C$ is regular on $V(\Gamma)$ because $\Gamma$ is $G$-vertex-transitive.
In this case, $|C||N|=|K|=|CK_v|=|C||K_v|$ and $|K_v|=|N|=|A_6|=2^3\cdot 3^2\cdot 5$, which is impossible by Proposition~\ref{prop=stabilizer}. Thus, $C_v\not=1$, and similarly, $N_v\not=1$. Since $\Gamma$ is $K$-symmetric
and has valency $5$, $K_v$ is primitive on the neighborhood $N(v)$, and since $C_v\unlhd K_v$
and $N_v\unlhd K_v$, both $C_v$ and $N_v$ are transitive on $N(v)$. It follows $5\di |N_v|$ and $5\di |C_v|$,
and hence $5^2\mid |K_v|$, contrary to Proposition~\ref{prop=stabilizer}. Now let $C\cap G\not=1$.
Then $G\leq C$ as $G$ is simple. Since $C\cong G$, we have $C=G$, and so $K=G\times N$. Set $L=C_X(N)$ and $M=LN$.
Then $L\unlhd X$, $G\leq L$ and $M=L\times N$ because $L\cap N=Z(N)=1$. Since $G\ntrianglelefteq X$,
$G$ is a proper subgroup of $L$ and hence $L_v\not=1$. It follows $5\di |L_v|$ because $\G$ is $X$-symmetric.
Since $|L||N|=|M|=|LM_v|=|L||M_v|/|L_v|$, we have $|M_v|=|L_v||N|$, which implies $5^2\mid |M_v|$,
contrary to Proposition~\ref{prop=stabilizer}.

The above contradictions imply $N\cap G\not=1$, and so $G\leq N$. Thus, $|N|/|G|=|N|/|N\cap G|\di 2^9\cdot 3^2\cdot 5$.
Since $G\ntrianglelefteq X$, $G$ is a proper subgroup of $N$, yielding $5\mid |N_v|$. It follows that $\Gamma$ is $N$-symmetric.
Recall that $N=T^s$ for some $s\geq 1$.

Suppose $s\geq 2$. We claim $|T|\di 2^9\cdot3^2\cdot5$. It suffices to show $|T|\di |N|/|G|$.
Note that $T\cap G=G$ or $1$. If $T\cap G=G$ then $G\leq T$ and so $|G|\di |T|$. Since $s\geq 2$,
we have $|T|\di |N|/|G|$. If $G\cap T=1$ then $|TG|=|T||G|$. Since $TG\leq N$, we have $|G|\di |T|^{s-1}$
and hence $|T|\di |N|/|G|$, as required. The proof also implies that if $s=2$ then $|G|\di |T|$.

Since $G\leq N$, the above claim means that both $G$ and $T$ are non-abelian simple $\{2,3,5\}$-groups.
Since $G\not=A_5$, Proposition~\ref{prop=235simplegroup} implies $G=A_6$ or $\PSU(4,2)$ and also $T=A_6$ or
$\PSU(4,2)$. Since $5^2\nmid |N|/|G|$, we have $s=2$ and so $|G|\di |T|$. Note that $N=T\times T$.
If both direct factors $T$ of $N$ are not semiregular, then $5\di |T_v|$ for each $T$ because $\G$ is $N$-symmetric,
yielding $5^2\di |N_v|$, contrary to Proposition~\ref{prop=stabilizer}. Thus, $N$ has a direct factor $T$ with $T_v=1$.
Since $|G|\di |T|$ and $\Gamma$ is $G$-vertex-transitive, $T$ is regular on $V(\Gamma)$ and so $|N_v|=|T|=|A_6|$ or $|\PSU(4,2)|$, that is, $2^3\cdot 3^2\cdot 5$ or $2^6\cdot 3^4\cdot 5$, of which both are impossible by Proposition~\ref{prop=stabilizer}.

The above contradictions imply $s=1$, that is, $N=T$ is a non-abelian simple group. Thus, $G\leq T\unlhd X$
and $|T:G|=|N|/|G|\di 2^9\cdot 3^2\cdot 5$. Since $\G$ is $X$-symmetric, it is also $T$-symmetric, that is, $5\di |T_v|$. Furthermore, $|T|/|G|=|T_v|/|G_v|$ and $G\not=T$ as $G\ntrianglelefteq X$. If $|T|/|G|\di 2^9\cdot3^2$, then $5\di |T_v|$, implying that $\G$ is $G$-symmetric. By Proposition~\ref{prop=cosetgraph}, $\G=\Cos(T,T_v,T_vtT_v)$ for some feasible $t$, and if $\G$ is $G$-symmetric then $\G=\Cos(G,G_v,G_vgG_v)$ for some feasible $g$ .

By the Frattini argument, $T=GT_v$. Let $H$ and $D$ be maximal subgroups of $T$ containing $G$ and $T_v$, respectively.
Then $T=HD$ is a maximal factorization of $T$. By Proposition~\ref{prop=simplegroupinedx2a3b5c}, $(G,T)$ is listed in Table~\ref{table=2}. Clearly, $(G,T)\neq (\rm PSL(2,7)$, $\rm PSU(3,3))$ as $5\nmid |\rm PSU(3,3)|$.
To finish the proof, we only need to show $(G,T)\not=(A_6, \PSU(4,2))$, $(A_7$, $\PSp(6,2))$,
$(A_7$, $A_9)$, $(A_8$, $\PSp(6,2))$, $(\PSL(2,11)$, $M_{11})$, $(\PSL(2,11)$, $M_{12})$, $(M_{11}$, $M_{12})$,
$(M_{23}$, $M_{24})$, $(\PSL(2,25)$, $^2F_{4}(2)')$, $(\PSL(2,16)$, $\PSp(4,4))$, $(\PSL(2,7)$, $\PSL(3,4))$,
$(\PSL(2,7)$, $A_7)$, $(\PSL(3,2)$, $A_8)$, $(\PSL(2,8)$, $\PSp(6,2))$, $(A_8$, $A_{10})$, $(A_7$, $A_{10})$,
$(A_9$, $\rm P\Omega^+(8,2))$ or $(\rm PSp(6,2)$, $\rm P\Omega^+(8,2))$.

Suppose $(G,T)=(A_6, \PSU(4,2))$, $(A_7, \PSp(6,2))$, $(A_7, A_9)$ or $(A_8, \PSp(6,2))$. Then $|T_v|/|G_v|=|T|/|G|\di 2^6\cdot3^2$ and $\G$ is $G$-symmetric. By Lemma~\ref{lem=S_9}, $G\unlhd A$, a contradiction.

Suppose $(G,T)=(\PSL(2,11),M_{11})$ or $(\PSL(2,11),M_{12})$. Then $|T_v|/|G_v|=|T|/|G|\di 2^4\cdot3^2$
and $\G$ is $G$-symmetric. By the Atlas \cite[p. 7]{Atlas} and Proposition \ref{prop=stabilizer},
$G_v\cong \mz_5$, $D_5$ or $A_5$. If $G_v\cong A_5$ then $|V(\G)|=11$, which is impossible because
there is no pentavalent graph of odd order. It follows $G_v\cong \mz_5$ or $D_5$,
and hence $|V(\G)|=132$ or $66$. By \cite[Theorem 4.1]{Guo12p} and \cite[Theorem 4.2]{hua}, $\G= \mathcal{G}_{132}^1$, $\mathcal{G}_{132}^2$, $\mathcal{G}_{132}^3$ or $\mathcal{G}_{66}$: $A=\PSL(2,11)\times \mz_2$
for $\G= \mathcal{G}_{132}^1$ and $A=\PGL(2,11)$ for other cases. In all cases, $G\unlhd A$, a contradiction.

Suppose $(G,T)=(M_{11},M_{12})$. Then $|T_v|/|G_v|=|T|/|G|\di 2^2\cdot3$ and $\G$ is $G$-symmetric.
By the Atlas \cite[p. 18]{Atlas} and Proposition \ref{prop=stabilizer}, $G_v\cong \mz_5$, $D_5$, $F_{20}$, $A_5$ or $S_5$,
and $|T_v|=|G_v||T|/|G|= 2^2\cdot3 \cdot 5$, $2^3\cdot3 \cdot 5$, $2^4\cdot 3\cdot 5$, $2^4\cdot 3^2\cdot 5$ or
$2^5\cdot 3^2\cdot 5$, respectively. By Proposition \ref{prop=stabilizer}, $(G_v,T_v)\cong (\mz_5,A_5)$, $(D_5,S_5)$,
$(A_5,A_4 \times A_5)$ or $(S_5,(A_4\times A_5)\rtimes \mz_2)$. By MAGMA, $T$ has no subgroups isomorphic to
$A_4 \times A_5$ or $(A_4 \times A_5)\rtimes \mz_2$, and $G$ has one conjugacy class of subgroups isomorphic
to $\mz_5$ or $D_5$. Take a given $G_v$ in each class: for $G_v\cong D_5$, there is no feasible $g$ for $G_v$,
and for $G_v\cong \mz_5$, there are $40$ feasible $g$ for $G_v$. For every such feasible $g$, the corresponding
coset graph $\G$ is a connected pentavalent $G$-symmetric graph with $A=M_{11}$. Thus, $G\unlhd A$, a contradiction.

Suppose $(G,T)=(M_{23},M_{24})$. Then $|T_v|/|G_v|=|T|/|G|\di 2^3\cdot3$ and $\G$ is $G$-symmetric.
By the Atlas \cite[p. 71]{Atlas} and Proposition \ref{prop=stabilizer}, $G_v\cong \mz_5$, $D_5$, $F_{20}$,
$A_5$, $S_5$, $\rm ASL(2,4)$, $\rm AGL(2,4)$, $\rm A \Sigma L(2,4)$ or $\rm A\Gamma L(2,4)$, and
$|T_v|=|G_v||T|/|G|= 2^3\cdot3 \cdot 5$, $2^4\cdot3 \cdot 5$, $2^5\cdot 3\cdot 5$, $2^5\cdot 3^2\cdot 5$,
$2^6\cdot 3^2\cdot 5$, $2^9\cdot 3^2\cdot 5$,  $2^9\cdot 3^3\cdot 5$, $2^{10}\cdot 3^2\cdot 5$ or
$2^{10}\cdot 3^3\cdot 5$, respectively. By Proposition~\ref{prop=stabilizer}, $(G_v,T_v)\cong (\mz_5,S_5)$,
$(A_5,(A_4 \times A_5)\rtimes \mz_2)$, $(S_5,\rm AGL(2,4))$, $(S_5,S_4 \times S_5)$ or
$(\ASL(2,4),\mathbb{Z}^6_2\rtimes  {\rm \Gamma L}(2,4))$. By MAGMA, $T$ has no subgroup isomorphic to
$S_4 \times S_5$, $\rm AGL(2,4)$ has no subgroup isomorphic to $S_5$, and
$\mathbb{Z}^6_2\rtimes  {\rm \Gamma L}(2,4)$ has no subgroup isomorphic to $\ASL(2,4)$ (both are subgroups of $M_{24}$).
Furthermore, $T$ has five conjugacy classes of subgroups isomorphic to $S_5$ and one conjugacy class of subgroups isomorphic to $(A_4\times A_5)\rtimes \mz_2$. Take a given  $T_v$ in each class, and computation shows that $T$ has no feasible $t$ for each $T_v$, a contradiction.

Suppose $(G,T)=(\PSL(2,25), {}^2F_{4}(2)')$. Then $|T_v|/|G_v|=|T|/|G|\di 2^8\cdot3^2$ and $\G$ is $G$-symmetric.
By the Atlas~\cite[p. 16]{Atlas} and Proposition~\ref{prop=stabilizer}, $G_v\cong \mz_5$, $D_5$, $F_{20}$, $A_5$
or $S_5$, and $|T_v|=|G_v||T|/|G|=2^8\cdot3^2 \cdot 5$, $2^9\cdot3^2 \cdot 5$, $2^{10}\cdot 3^2\cdot 5$,
$2^{10}\cdot 3^3\cdot 5$ or $2^{11}\cdot 3^3\cdot 5$, respectively. By Proposition~\ref{prop=stabilizer},
$(G_v,T_v)\cong (D_5,\mathbb{Z}^6_2\rtimes  {\rm \Gamma L}(2,4))$, and by MAGMA, $G$ contains two conjugacy
classes of subgroups isomorphic to $D_5$. By taking a given $T_v$ in each conjugacy class, computation shows
that for each $G_v$ there are $40$ feasible $g$, and for every such $g$, the corresponding coset graph $\G$
is a connected pentavalent $G$-symmetric graph with $A=\PSL(2,25)\times \mz_2$. Thus, $G\unlhd A$, a contradiction.

Suppose $(G,T)=(\PSL(2,16),\PSp(4,4))$. Then $|T_v|/|G_v|=|T|/|G|=2^4 \cdot 3\cdot5\di |T_v|$.
Since $|T|=2^8 \cdot 3^2 \cdot5^2 \cdot17$, Proposition~\ref{prop=stabilizer} implies $T_v\cong A_4 \times A_5$,
$(A_4 \times A_5)\rtimes \mz_2$, $S_4 \times S_5$, $\ASL(2,4)$, $\rm AGL(2,4)$, $\rm A\Sigma L(2,4)$, $\rm A\Gamma L(2,4)$.
However, by MAGMA, $\PSp(4,4)$ has no subgroup isomorphic to $(A_4 \times A_5)\rtimes \mz_2$, $S_4 \times S_5$, $\ASL(2,4)$,
$\rm AGL(2,4)$, $\rm A\Sigma L(2,4)$ or $\rm A\Gamma L(2,4)$. For $T_v\cong  A_4 \times A_5$, $T_v$ has two conjugacy classes in $T$. By taking a given $T_v$ in each conjugacy class, computation shows that
there is no feasible $t$, a contradiction.

Suppose $(G,T)=(\PSL(2,7),\PSL(3,4))$. Then $|T|=2^6 \cdot 3^2 \cdot5 \cdot7$ and $|T_v|/|G_v|=|T|/|G|=2^3 \cdot 3\cdot5\di |T_v|$. Note that
there is no pentavalent graph of odd order.
By Proposition~\ref{prop=stabilizer}, $T_v\cong S_5$, $A_4 \times A_5$ or $(A_4 \times A_5)\rtimes \mz_2$.
By the Atlas \cite[p. 23]{Atlas}, $T$ has no subgroup isomorphic to $S_5$, $A_4 \times A_5$
or $(A_4 \times A_5)\rtimes \mz_2$, a contradiction.

Suppose $(G,T)=(\PSL(2,7),A_7)$. Then $|T_v|/|G_v|=|T|/|G|= 3\cdot5\di |T_v|$. By the Atlas \cite[p. 10]{Atlas}
and Proposition~\ref{prop=stabilizer}, $T_v\cong A_5 $. By MAGMA, $T_v \cong A_{5}$ has two conjugacy classes in $T$
and by taking a given $T_v$ in each conjugacy class, computation shows that there is no feasible $t$, a contradiction.

Suppose $(G,T)=(\PSL(3,2),A_8)$. Then $|T_v|/|G_v|=|T|/|G|= 2^3 \cdot3\cdot5\di |T_v|$. By the Atlas \cite[p. 22]{Atlas}
and Proposition~\ref{prop=stabilizer}, $T_v\cong S_5 $. By MAGMA, $T_v \cong S_5$ has two conjugacy classes and
by taking a given $T_v$ in each conjugate class, computation shows that there is no feasible $t$, a contradiction.

Suppose $(G,T)=(\PSL(2,8),\PSp(6,2))$. Then  $|T_v|/|G_v|=|T|/|G|=2^6 \cdot 3^2\cdot5\di |T_v|$.
Since $|T|=2^9 \cdot 3^4 \cdot5 \cdot7$, Proposition~\ref{prop=stabilizer} implies  $T_v\cong S_4 \times S_5$,
$\rm AGL(2,4)$ or $\rm A\Gamma L(2,4)$. However, by the Atlas~\cite[p. 46]{Atlas}, $\PSp(6,2)$ has no subgroup
isomorphic to $S_4 \times S_5$, $\rm AGL(2,4)$ or $\rm A\Gamma L(2,4)$, a contradiction.

Suppose $(G,T)=(A_8,A_{10})$ or $(A_7,A_{10})$. Then $3^2 \di |T|/|G|=|T_v|/|G_v|$, and by the Atlas \cite[p. 48]{Atlas}
and Proposition~\ref{prop=stabilizer}, $T_v\cong A_4 \times A_5$ or $(A_4 \times A_5)\rtimes \mz_2$. For each case,
by MAGMA, $T_v$ has two conjugacy classes. By taking a given $T_v$ in each conjugacy class, computation shows that
there is no feasible $t$, a contradiction.

Suppose $(G,T)=(A_9,\rm P\Omega^+(8,2))$ or $(\rm PSp(6,2),\rm P\Omega^+(8,2))$. Then $|T_v|/|G_v|=|T|/|G|=2^6\cdot 3\cdot 5$
or $2^3\cdot 3\cdot 5$. By Proposition~\ref{prop=stabilizer}, $T_v\cong S_5$, $A_4\times A_5 $, $(A_4\times A_5)\rtimes \mz_2$,
$S_4 \times S_5$, $\rm ASL(2,4)$, $\rm AGL(2,4)$, $\rm A\Sigma \rm L(2,4)$, $\rm A \Gamma \rm L(2,4)$ or $\mz_2^6\rtimes \Gamma \rm L(2,4)$. By MAGMA, $\rm P\Omega^+(8,2)$ has no subgroup isomorphic to $S_4 \times S_5$, $\rm ASL(2,4)$, $\rm AGL(2,4)$, $\rm A\Sigma \rm L(2,4)$, $\rm A \Gamma \rm L(2,4)$ or $\mz_2^6\rtimes \Gamma \rm L(2,4)$. Thus, $T_v\cong S_5$, $A_4\times A_5 $ or $(A_4\times A_5)\rtimes \mz_2$, and the computation shows that $T_v$ has $15$, $3$, $3$ conjugacy classes, respectively. By taking a given $T_v$ in each conjugacy class, there is no feasible $t$ respectively, a contradiction.
\hfill\qed

\begin{lem}\label{lem=T normal}
Let $\G$ be a connected pentavalent symmetric $G$-vertex-transitive graph and
let $\Aut(\G)$ have a non-trivial radical $R$ with at least three orbits. Assume $G\ntrianglelefteq \Aut(\G)$
and $RG=R\times G$. Then $\Aut(\G)$ contains a non-abelian simple normal subgroup $T$ such that $G<T$ and $(G,T)=(\Omega^-_8(2),\PSp(8,2))$, $(A_{14},A_{16})$, $(\PSL(2,8),A_{9})$ or $(A_{n-1},A_n)$ with $n\geq 8$ and $n\di 2^9\cdot3^2\cdot5$.
\end{lem}

\f {\bf Proof:} Set $A=\Aut(\G)$ and $B=RG=R\times G$. Since $R\neq 1$ has at least three orbits,
by Proposition~\ref{prop=atlesst3orbits}, the quotient graph $\G_{R}$ is a connected pentavalent
$A/R$-symmetric graph with $A/R\leq \Aut(\G_{R})$. Furthermore, $\G_R$ is $B/R$-vertex-transitive and $G$
is characteristic in $B$. It follows that $B\ntrianglelefteq A$ because $G\ntrianglelefteq A$, and
hence $G\cong B/R\ntrianglelefteq A/R$. Since $R$ is the largest solvable normal subgroup of $A$,
$A/R$ has trivial radical. By Lemma~\ref{lem=insolvable2}, $A/R$ has a non-abelian simple normal subgroup $I/R$ such that $B/R\leq I/R\cong T$ and $(B/R,I/R)\cong (G,T)=(\Omega^-_8(2),\PSp(8,2))$, $(A_{14},A_{16})$, $(\PSL(2,8),A_{9})$ or $(A_{n-1},A_n)$ with $n\geq 6$ and $n\di 2^9\cdot3^2\cdot5$. If $(G,T)\not=(A_5,A_6)$ then $A=S_6$, contrary to the fact that $A$ has a non-trivial radical. Thus, $n\geq 8$

Let $C=C_I(R)$. Then $G\leq C$ as $B=R\times G\leq I$. Clearly, $C\unlhd I$ and $C\cap R\leq Z(C)$.
Since $G\cong R\times G/R\leq CR/R$ and $C/C\cap R\cong CR/R\unlhd I/R\cong T$, we have  $C\cap R=Z(C)$, $C/Z(C)\cong T$
and $I=CR$. Thus, $C'/C'\cap Z(C)\cong C'Z(C)/Z(C)=(C/Z(C))'=C/Z(C)\cong  T$, and so  $Z(C')=C'\cap Z(C)$, $C=C'Z(C)$ and $C'/Z(C')\cong T$. Furthermore, $C'=(C'Z(C))'=C''$, implying that $C'$ is a covering group of $T$. It follows $G\leq C'$ because $C/C'$ is abelian.

By \cite[Theorem 5.1.4]{Kleidman} and the Atlas~\cite[p. 123]{Atlas}, $\Mult(A_n)=\mz_2$ for $n\geq 8$ and $\Mult(\PSp(8,2))=1$. Since $C'$ is a covering group of $T$, we have $Z(C')=1$ or $\mz_2$.

Suppose $Z(C')=\mz_2$. Then $(G,T)=(A_{14},A_{16})$, $(\PSL(2,8),A_9)$, or $(A_{n-1},A_n)$ with $n\geq 8$
and $n\di 2^9\cdot 3^2\cdot 5$. Furthermore, $C'$ is the unique double cover of $A_n$ with $n\geq 8$,
that is, $C'=2.A_n$. Since $Z(C')\unlhd I$ and $I/R\cong T$ is simple, we have $Z(C')\leq R$, and hence
$G\times Z(C')\leq C'$ as $B=G\times R$.

Let $(G,T)=(A_{14},A_{16})$. Then $C'=2.A_{16}$. Clearly, $A_{14}$ has no proper subgroup of index less than $13$, and by MAGMA, it also has no subgroup of index $15$ or $16$. This implies that each subgroup of $A_{16}$ isomorphic to $A_{14}$ fixes exactly two vertices in $\{1,2,\cdots, 16\}$, and so lies in a subgroup of index $16$ of  $A_{16}$. It follows that $G\times Z(C')$ lies in a subgroup of index $16$ of  $C'$, that is, $G\times Z(C')\leq L$ and $|C':L|=16$ for some $L\leq C'$. By Proposition~\ref{prop=covering group}, $L\cong 2.A_{15}$ and $G\times Z(C')\cong \mz_2\times A_{14}$, which is impossible because  $G\times Z(C')\cong A_{14}\times\mz_2$.

Let $(G,T)=(\PSL(2,8),A_9)$. Then $C'=2.A_{9}=G\cdot C'_v$, and so $3^2 \di |C'_v|$ as $|C'|=2^7\cdot 3^4\cdot 5\cdot7$. By Proposition~\ref{prop=stabilizer}, $C_v'\cong A_4\times A_5$,
$(A_4\times A_5)\rtimes \mz_2$, $S_4\times S_5$, $\rm ASL(2,4)$ or $\rm AGL(2,4)$. Note that $C'$ has a subgroup $H$ isomorphic to $2.A_6$ with $Z(C')\leq H$. Since $H\cong \rm SL(2,9)$, $Z(C')$ is the unique subgroup of order $2$ in $H$. Since $Z(C')\unlhd C'$, we have $Z(C')\nleq C'_v$. Suppose that $|H_w|$ is odd for any $w\in V(\G)$.
Then $|w^H|=|H|/|H_w|=2^4\cdot (3^2\cdot 5/|H_w|)$, implying that the orbit of $H$ containing $w$ has length
a multiple of $2^4$, and so $|V(\G)|$ is a multiple of $2^4$, contrary to the fact
$|V(\G)|\di|\PSL(2,8)|=2^3\cdot 3^2\cdot 7$. Thus, we may assume that $H_v$ is even.
Since $H_v=H\cap C'_v$ and  $Z(C')$ is the unique subgroup of order $2$ in $H$, we have
$Z(C')\leq H_v$ and hence $Z(C')\leq C'_v$, a contradiction.

Let $(G,T)=(A_{n-1},A_n)$ with $n\geq 8$. Then $G\times Z(C')$ is a subgroup of index $n$ of $C'$ isomorphic to $A_{n-1}\times \mz_2$, which is impossible by Proposition~\ref{prop=covering group}.

The above contradictions give rise to $Z(C')=1$. Then $C'\cong T$ and $G<C'\unlhd I$.
Since $|I|=|I/R||R|=|T||R|=|C'||R|$ and $C'\cap R=1$, we have $I=C'\times R$. Then $C'$ is characteristic in $I$,
and so $C'\unlhd A$ because $I\unlhd A$. It follows that $A$ has a non-abelian simple normal subgroup $C'$ such that $G< C'$ and $(G,C')\cong (G,T)=(\Omega^-_8(2),\PSp(8,2))$, $(A_{14},A_{16})$, $(\PSL(2,8),A_{9})$ or $(A_{n-1},A_n)$ with $n\geq 8$ and $n\di 2^9\cdot3^2\cdot5$.
\hfill\qed

Now, we are ready to prove Theorem \ref{theo=main}.

\medskip

\f {\bf The proof of Theorem \ref{theo=main}:} Let $G$ be a non-abelian simple group and $\G$ a connected pentavalent
symmetric $G$-vertex-transitive graph with $v\in V(\G)$. Let $A=\Aut(\G)$ and $R$ the radical of $A$. To prove the theorem,
we may assume $G\ntrianglelefteq A$. If $R=1$, the theorem is true by  Lemma~\ref{lem=insolvable2}. In what follows
we assume $R\not=1$. If $G=A_5$, by Lemma~\ref{lem=A_5}, $\G=K_6$ and $A=S_6$, which is impossible because $A$ has trivial radical. Thus, $G\not=A_5$.

Set $B=RG$. Then $G\cap R=1$ and $|B|=|R||G|$. By Lemma~\ref{lem=GH},
$|R|=|B_v|/|G_v|\di 2^9\cdot 3^2\cdot 5$.

Suppose that $R$ has one or two orbits on $V(\G)$. Since $|G|=|V(\G)||G_v|$,
we have $|R|=|v^R||R_v|=|G||R_v|/|G_v|$ or $|G||R_v|/(2\cdot|G_v|)$, and since $R$ and $G_v$ are $\{2,3,5\}$-groups,
$G$ is a non-abelian $\{2,3,5\}$-simple group. By Proposition~\ref{prop=235simplegroup}, $G=A_6$ or $\rm PSU(4,2)$ as $G\neq A_5$, and since $|R|=|B_v|/|G_v|$, we have $|B_v|=|G||R_v|$ or $|G||R_v|/2$. Note that $5\di |G|$, and so $5\di |B_v|$, that is, $\G$ is $B$-symmetric.
If $R_v=1$, then $|B_v|=|G|$ or $|G|/2$, that is, $|B_v|=2^3\cdot 3^2 \cdot 5$, $2^2\cdot 3^2 \cdot 5$,
$2^6\cdot 3^4 \cdot 5$ or $2^5\cdot 3^4 \cdot 5$, of which all are impossible by Proposition~\ref{prop=stabilizer}.
If $R_v\not=1$ then $5\di |R_v|$ as $R\unlhd A$, and in this case, $5^2\di |B_v|=|G||R_v|$, a contradiction.

Now we have shown that $R$ has at least three orbits. By Lemma~\ref{lem=T normal}, if $B=R\times G$ then the theorem is true. To finish the proof, we only need to show $B=R\times G$.

By Proposition~\ref{prop=atlesst3orbits}, $R$ is semiregular on $V(\G)$ and the quotient graph $\G_{R}$ is a
connected pentavalent $A/R$-symmetric graph with $A/R\leq \Aut(\G_{R})$. Moreover, $\G_R$ is $B/R$-vertex-transitive.
Since $R$ is the largest solvable normal subgroup of $A$, $A/R$ has trivial radical, by Lemma~\ref{lem=insolvable2},
$B/R\unlhd A/R$ or $A/R$ has a non-abelian simple normal subgroup $I/R$ such that $B/R\leq I/R$ and
$(B/R,I/R) \cong (G,T)$ with $(G,T)=(\Omega^-_8(2),\PSp(8,2))$, $(A_{14},A_{16})$, $(\PSL(2,8),A_{9})$ or $(A_{n-1},A_n)$ with $n\geq 6$ and $n\di 2^9\cdot3^2\cdot5$.

Suppose $B/R \unlhd A/R$. Then $B\unlhd A$ and $5\di |B_v|$. Note that $1\not=|R|=|B_v|/|G_v|$. Then $5 \di |R|$
or $5\di |G_v|$. By Lemmas~\ref{lem=RG} and \ref{lem=GP}, $B=R\times G$, that is, $G$ is characteristic in $B$, and hence $G\unlhd A$,
contrary to the assumption $G\ntrianglelefteq A$. Thus, $B/R\leq I/R\unlhd A/R$. It follows $(I/R)_\a\not=1$ for $\a\in V(\G_R)$ and hence $5\di |(I/R)_\a|$. Thus, $\G_R$ is $I/R$-symmetric. Since $|V(\G_R)|=|V(\G)|/|R|=|G|/(|R||G_v|)$ and $|T|=|I/R|=|V(\G_R)||(I/R)_\a|=|G|/(|R||G_v|)\cdot |(I/R)_\a|$, we have $|(I/R)_\a|=|R||T||G_v|/|G|$, and by  Proposition~\ref{prop=stabilizer}, $|R||T||G_v|/|G|$ is a divisor of $2^9\cdot 3^2\cdot 5$.

Since $|R|\di 2^9\cdot 3^2\cdot 5$, we may write $|R|=2^m\cdot3^n\cdot 5^k$, where $0\leq m\leq 9$,
$0\leq n\leq 2$ and $0\leq k\leq 1$. Since $R$ is solvable, there exists a series of subgroups of $B$:
\begin{center}
$B>R=R_{s}>\cdots >R_1>R_0=1$
\end{center}
such that $R_i\unlhd B$ and $R_{i+1}/R_i$ is an elementary abelian $r$-group with $0\leq i\leq s-1$, where
$r=2$, $3$ or $5$. Clearly, $G\leq B$ has a natural action on $R_{i+1}/R_i$ by conjugation.

Suppose to the contrary $B\neq R\times G$. Then there exists some $0\leq j\leq s-1$ such that $GR_i=G\times R_i$ for any $0\leq i\leq j$,
but $GR_{j+1}\not=G\times R_{j+1}$. If $G$ acts trivially on $R_{j+1}/R_j$ by conjugation, then $[GR_j/R_j,R_{j+1}/R_j]=1$.
Since $GR_j/R_j\cong G$ is simple, we have $(GR_j/R_j)\cap (R_{j+1}/R_j)=1$. Note that $|GR_{j+1}/R_j|=|GR_{j+1}/R_{j+1}||R_{j+1}/R_j|=|G||R_{j+1}/R_j|=|GR_j/R_j||R_{j+1}/R_j|$. Then $GR_{j+1}/R_j=GR_j/R_j\times R_{j+1}/R_j$. In particular, $GR_j\unlhd GR_{j+1}$ and so $G\unlhd GR_{j+1}$ because $GR_j=G\times R_j$
implies that $G$ is characteristic in $GR_j$. It follows that $GR_{j+1}=G\times R_{j+1}$, a contradiction.
Thus, $G$ acts non-trivially on $R_{j+1}/R_j$, and the simplicity of $G$ implies that $G$ acts faithfully
on $R_{j+1}/R_j$.

Set $|R_{j+1}/R_j|=r^\ell$. Since $|R|=2^m\cdot3^n\cdot 5^k$, we have $\ell\leq m\leq 9$ for $r=2$, $\ell\leq 2$ for $r=3$
and $\ell\leq 1$ for $r=5$. Recall that $G$ acts faithfully on  $R_{j+1}/R_j$ and $R_{j+1}/R_j$ is elementary abelian.
Then $G\leq \GL(\ell,r)$. If $r=3$ or $5$ then $G\leq \GL(2,3)$ or $\GL(1,5)$, yielding that $G$ is solvable, a contradiction.
Thus, $G\leq \GL(m,2)=\PSL(m,2)$ with $m\leq 9$. Recall that $(G,T)=(\Omega^-_8(2),\PSp(8,2))$, $(A_{14},A_{16})$, $(\PSL(2,8),A_{9})$ or $(A_{n-1},A_n)$ with $n\geq 6$ and $n\di 2^9\cdot3^2\cdot5$.

Let $G=A_{14}$ or $A_{n-1}$ with $n\di 2^9\cdot3^2\cdot5$ and $n\geq 12$. Since $11\di |A_{n-1}|$ and $11\nmid |\PSL(9,2)|$, we have $G\nleq \PSL(9,2)$, a contradiction.

Let $G=\PSL(2,8)$ or $A_7$, since $|G|=|V(\G)||G_v|$ and $R$ is semiregular on $V(\G)$, we have $|R|\di |G|/|G_v|$, and so $m\leq 2$ because $\G_R$ cannot be of odd order. It follows that $A_7\leq \PSL(2,2)$ or $\PSL(2,8)\leq \PSL(2,2)$, a contradiction. Similarly, for $G=A_9$, we have $m\leq 5$ and $A_9\leq \PSL(5,2)$, which is also impossible because $3^4\di |A_9|$ and $3^4\nmid |\PSL(5,2)|$.

Since $G\not=A_5$, the left pairs are $(G,T)=(A_8,A_9)$ or $(\Omega^-_8(2)$, $\PSp(8,2))$.
For $(G,T)=(A_8,A_9)$, since $|(I/R)_\a|=|R||T||G_v|/|G|$ and $5\di |(I/R)_\a|$, $|R|$ or $|G_v|$
is divisible by $5$, by Lemmas~\ref{lem=RG} and \ref{lem=GP}, $B=R\times G$, a contradiction. Thus,
$(G,T)=(\Omega^-_8(2),\PSp(8,2))$. Since $|T|/|G|=2^4\cdot3\cdot5$ and $|R||T||G_v|/|G|$ is a divisor
of $2^9\cdot 3^2\cdot 5$, $|R|$ is a divisor of $2^5\cdot3$. Then $m\leq 5$ and $G\leq \PSL(5,2)$,
which is impossible because $17\di |\Omega^-_8(2)|$ and $17\nmid |\PSL(5,2)|$.

The above contradictions imply $B=R\times G$, as required. This completes the proof of Theorem~\ref{theo=main}.

%
%

\hfill\qed

\f {\bf The proof of Corollary \ref{cor=arc}:} Let $G$ be a non-abelian simple group
and $\G$ a connected pentavalent $G$-symmetric graph with $v\in V(\G)$. Let $A=\Aut(\G)$ and let $R$ the radical of $A$.
To prove the corollary, we may assume $G\ntrianglelefteq A$. By Theorem~\ref{theo=main}, $A$ has a non-abelian simple normal subgroup $T$
such that $G<T$ and $(G,T)$ is as given in Theorem~\ref{theo=main}. Since $G$ is symmetric and $G< T\unlhd A$,
$ |G_v|$ and $ |T_v|$ are divisible by $5$, and since $|T|/|G|=|T_v|/|G_v|$, by Proposition~\ref{prop=stabilizer} $(G,T)=(A_{n-1},A_n)$ with $n\geq 6$, $n\di 2^9\cdot3^2$ and $n\neq 2^9,2^8,2^7,2^6,2^9\cdot3,2^8 \cdot3$.

Let $(G,T)=(A_{2^3-1},A_{2^3})$ or $(A_{3^2-1},A_{3^2})$. By Lemma~\ref{lem=S_9}, $G\unlhd A$, a contradiction.

Let $(G,T)=(A_{2\cdot 3^2-1},A_{2\cdot 3^2})$. By Proposition~\ref{prop=stabilizer},
$(G_v,T_v)\cong (F_{20}\times \mz_2, A_4 \times A_5)$ or $(F_{20}\times \mz_4, (A_4 \times A_5)\rtimes \mz_2$).
However, by MAGMA, $A_4\times A_5$ or $(A_4\times A_5)\rtimes \mz_2$ has no subgroup isomorphic to $F_{20}\times \mz_2$ or $F_{20}\times \mz_4$, respectively, a contradiction.

Let $T=A_6$. Then $G=A_5$ and by Lemma~\ref{lem=S_9}, $\Gamma=K_6$ and $A=S_6$. To finish the proof,
we only need to show that for $n=2^5, 2^7 \cdot 3, 2^7 \cdot 3^2, 2^8 \cdot 3^2$ and $2^9\cdot 3^2$,
we have $A=T$. To do this, it suffices to show $A_v=T_v$ because $A=TA_v$.

Let $n=2^9\cdot 3^2$ or $2^8 \cdot 3^2$. By  Proposition~\ref{prop=stabilizer}, $(G_v,T_v)\cong (\mz_5,\mz_2^6\rtimes \rm \Gamma L(2,4))$ or $(D_5,\mz_2^6\rtimes \rm \Gamma L(2,4))$, respectively. This forces $A_v=T_v$ because $T_v\leq A_v$, as required.

Let $n=2^7\cdot 3^2$. By Proposition~\ref{prop=stabilizer}, $(G_v,T_v)\cong (\mz_5,\rm A \Gamma L(2,4))$,
$(F_{20},\mz_2^6\rtimes \rm \Gamma L(2,4))$ or $(D_{10},\mz_2^6\rtimes \rm \Gamma L(2,4))$.
Clearly, if $(G_v,T_v)\cong (F_{20},\mz_2^6\rtimes \rm \Gamma L(2,4))$ or $(D_{10},\mz_2^6\rtimes \rm \Gamma L(2,4))$,
then $T_v=A_v$. Assume $(G_v,T_v)\cong (\mz_5,\rm A \Gamma L(2,4))$. If $T_v\not=A_v$,
then by Proposition~\ref{prop=stabilizer} we have $A_v=\mz_2^6\rtimes \rm \Gamma L(2,4)$,
and since $T\trianglelefteq A$, we have $T_v\unlhd A_v$. This implies that $\mz_2^6\rtimes \rm \Gamma L(2,4)$ has a normal subgroup isomorphic to $\rm A\Gamma L(2,4)$, which is impossible by~\cite[Theorem 1.1]{Weiss} and MAGMA. Thus, $T_v=A_v$, as required.

Let $n=2^7 \cdot 3$. By Proposition~\ref{prop=stabilizer}, $(G_v,T_v)\cong (\mz_5,\rm A \Sigma L(2,4))$ or $(A_5,\mz_2^6\rtimes \rm \Gamma L(2,4))$. For the latter, $A_v=T_v$. Assume $(G_v,T_v)\cong (\mz_5,\rm A \Sigma L(2,4))$. If $A_v\not=T_v$ then by Proposition~\ref{prop=stabilizer} we have $A_v=\mz_2^6\rtimes \rm \Gamma L(2,4)$ or $\rm A \Gamma L(2,4)$, which are impossible because $\mz_2^6\rtimes \rm \Gamma L(2,4)$ and $\rm A \Gamma L(2,4)$ have no normal subgroup isomorphic to $\rm A\Sigma L(2,4)$ by~\cite[Theorem 1.1]{Weiss} and MAGMA. Thus, $T_v=A_v$, as required.

Let $n=2^5$. By Proposition \ref{prop=stabilizer}, $(G_v,T_v)\cong (A_5,\rm A \Sigma L(2,4))$ or $(A_4 \times A_5,\mz_2^6\rtimes \rm \Gamma L(2,4))$. A similar argument to the above paragraph implies $T_v=A_v$, as required. This completes the proof of Corollary~\ref{cor=arc}.
\hfill\qed

\f {\bf The proof of Corollary \ref{cor=regular}:} Let $G$ be a non-abelian simple group and $\G$ a connected pentavalent $G$-regular graph with $v\in V(\G)$. Let $A=\Aut(\G)$ and $R$ the radical of $A$. To prove the corollary, we may assume $G\ntrianglelefteq A$. By Theorem~\ref{theo=main}, $A$ contains non-abelian simple normal subgroup $T$ containing $G$ such that $(G,T)$ is given in Theorem \ref{theo=main}. Since $G$ is regular and $G< T\unlhd A$, $G_v=1$ and $T$ is symmetric. Since $|T|/|G|=|T_v|/|G_v|=|T_v|$, by Proposition~\ref{prop=stabilizer}, $(G,T)=(A_{n-1},A_n)$ with $n=2\cdot 5$, $2^2\cdot 5$, $2^3\cdot5$, $2\cdot3\cdot 5$, $2^4\cdot5$, $2^3\cdot 3\cdot5$, $2^4\cdot 3^2\cdot5$, $2^6\cdot 3\cdot5$, $2^5\cdot 3^2\cdot5$, $2^7\cdot 3\cdot5$, $2^6\cdot 3^2\cdot5$, $2^7\cdot 3^2\cdot5$ or $2^9\cdot 3^2\cdot5$.
\hfill\qed
\medskip

To end the paper, we give an example to show that the pair $(G,T)=(A_{39},A_{40})$ in Corollary~\ref{cor=regular} can happen.

\begin{exam}\label{NonNormalExample} Let $G=A_{39}$ and $T=A_{40}$. Define $x,y,z,g\in T$ as following:

\begin{enumerate}
\item[]$x=(1$ $2$ $7$ $18$ $5)(3$ $11$ $8$ $17$ $14)(4$ $9$ $22$ $33$ $15)(6$ $20$ $19$ $10$ $21)$ $(12$ $28$ $26$ $31$ $23)$
    \\
   \mbox{}\hskip 0.8cm $(13$ $27$ $24$ $32$ $30)$$(16$ $35$ $34$ $25$ $36)$$(29$ $39$ $38$ $40$ $37)$,\\
$y=(1$ $3$ $12$ $6)$$(2$ $8$ $23$ $10)$$(4$ $13$ $29$ $16)(5$ $17$ $28$ $19)(7$ $14$ $31$ $20)(9$ $24$ $37$ $25)$\\
\mbox{}\hskip 0.8cm $(11$ $26$ $21$ $18)(15$ $32$ $39$ $34)$ $(22$ $30$ $40$ $35)(27$ $38$ $36$ $33)$,\\
$z=(1$ $4)$$(2$ $9)$$(3$ $13)$$(5$ $15)$$(6$ $16)$$(7$ $22)$$(8$ $24)$$(10$ $25)$$(11$ $27)$$(12$ $29)$$(14$ $30)$ $(17$ $32)$\\
\mbox{}\hskip 0.8cm $(18$ $33)$$(19$ $34)$$(20$ $35)$$(21$ $36)$$(23$ $37)$$(26$ $38)$$(28$ $39)$$(31$ $40)$,\\
$g=(1$ $7)$$(2$ $24)$$(3$ $14)$$(4$ $40)$$(5$ $11)$$(6$ $20)$$(8$ $37)$$(9$ $10)$$(12$ $31)$$(13$ $35)$$(15$ $36)$$(16$ $30)$ \\
\mbox{}\hskip 0.8cm $(17$ $26)$$(18$ $19)$$(21$ $28)$$(22$ $29)$$(23$ $25)$$(27$ $39)$$(32$ $33)$$(34$ $38)$.\\
\end{enumerate}

\rm {By MAGMA\cite{magma}, $H=\langle x,y,z\rangle\cong F_{20}\times \mz_2$, $T=\langle H,g\rangle$, $|H:H\cap H^g|=5$ and $H$ is regular on $\{1,2,\cdots, 40\}$. Thus, $T=GH$ with $G\cap H=1$ and the coset graph $\G=\Cos(T,H,HgH)$ is a connected pentavalent $T$-symmetric $G$-regular graph, where $G$ and $T$ are viewed as groups of automorphisms of $\G$ by right multiplication.

By Corollary~\ref{cor=regular}, $T\unlhd \Aut(\G)$. Again by MAGMA~\cite{magma}, $\Aut(T,H,HgH)\cong \tilde{H}$, where $\tilde{H}$ is the automorphism group of $T$ induced by conjugation of elements in $H$. Thus, $\Aut(\G)=T$ by \cite[Lemma 2.10]{Wang}.}
\end{exam}

\medskip
\f {\bf Acknowledgement:} This work was supported by the National Natural Science Foundation of China (11571035, 11231008, 11271012) and by the 111 Project of China (B16002).


\begin{thebibliography}{99}
\bibitem{magma}
W. Bosma, J. Cannon, and C. Playoust, The MAGMA algebra system I: The user language, J. Symbolic Comput. 24 (1997), 235-265.
\bibitem{Atlas}
J. H. Conway, R. T. Curtis, S. P. Norton, R. A. Parker, and R. A Wilson, Atlas of Finite Group,
Clarendon Press, Oxford, 1985.
\bibitem{FangLW}
X.G. Fang, L.J. Jia, J. Wang, On the automorphism groups of symmetric graphs admitting an almost simple group,
Europ. J. Combin. 29 (2008) 1467-1472.
\bibitem{Fang4}
X. G. Fang, C. H. Li, and M. Y. Xu, On edge-transitive Cayley graphs of valency four,
Europ. J. Combin. 25 (2004), 1107-1116.
\bibitem{Fang5}
X. G. Fang, X. S. Ma, and J. Wang, On locally primitive Cayley graphs of finite simple groups,
J. Combin. Theory Ser. A 118 (2011), 1039-1051.
\bibitem{FangP}
X. G. Fang and C. E. Praeger, Finite two-arc transitive graphs admitting a Suzuki simple group, Comm. Algebra 27 (1999), 3727-3754.
\bibitem{FP2}
X.G. Fang and C.E. Praeger, Finite two-arc-transitive graphs admitting a Ree simple group, Comm. Algebra 27 (1999) 3755-3769.
\bibitem{FangP2}
X. G. Fang, C. E. Praeger, and J. Wang, On the automorphism group of Cayley graphs of finite simple groups, J. London Math. Soc. 66 (2002), 563-578.
\bibitem{FengLZ}
Y.-Q. Feng, C.H. Li, J.-X. Zhou, Symmetric cubic graphs with sovable automorphism groups, Europ. J. Combin. 45 (2015) 1-11.
\bibitem{Godsil}
C. D. Godsil, On the full automorphism group of a graph,
Combinatorica 1 (1981), 243-256.
\bibitem{Guo}
S.-T. Guo and Y.-Q. Feng, A note on pentavalent $s$-transitive graphs,
Discrete Math. 312 (2012), 2214-2216.
\bibitem{Guo12p}
S.-T. Guo, J.-X. Zhou, and Y.-Q. Feng, Pentavalent symmetric graphs of order $12p$,
Electronic J. Combin. 18 (2011), 1-13 $\sharp$P233.
\bibitem{Hall}
P. Hall, A note on soluble groups,
J. London Math. Soc 3 (1928), 98-105.
\bibitem{hua}
X.-H. Hua, Y.-Q. Feng, and J. Lee, Pentavalent symmetric graphs of order $2pq$,
Discrete Math. 311 (2011), 2259-2267.
\bibitem{Huppert}
B. Huppert, Eudiche Gruppen I, Springer-Verlag, 1967.
\bibitem{HuppertB}
B. Huppert and N. Blackburn, Finite Groups III, Springer-Verlag, 1982.
\bibitem{Huppert2}
B. Huppert and W. Lempken, Simple groups of order divisible by at most four primes,
Proc. of the F. Scorina Gemel State University 16 (2000), 64-75.
\bibitem{Kleidman}
P. Kleidman and M. Liebeck, The Subgroup Structure of The Finite Classical Groups,
Cambridge Univ. Press, Cambridge, 1990.
\bibitem{CHLi}
C. H. Li, Isomorphisms of finite Cayley graphs,
Ph.D. Thesis, The University of Western Australia, 1996.
\bibitem{LiLZ}
C. H. Li, Z. P. Lu, H. Zhang, Tetravalent edge-transitive Cayley graphs with odd number of vertices,
J. Combin. Theory Ser. B 96 (2006)\textcolor[rgb]{1,0,0}{,} 164-181.
\bibitem{Lorimer}
P. Lorimer, Vertex-transitive graphs: symmetric graphs of prime valency,
J. Graph Theory 8 (1984), 55-68.
\bibitem{Miller}
R. C. Miller, The trivalent symmetric graphs of girth at most six,
J. Combin. Theory Ser. B 10 (1971), 163-182.
\bibitem{P1} C.E. Praeger, On a reduction theorem for finite, bipartite, $2$-arc-transitive graphs,
Australas. J. Combin. 7 (1993), 21-36.
\bibitem{P3} C.E. Praeger, An O'Nan-Scott theorem for finite quasiprimitive permutation groups and an application to $2$-arc transitive graphs, J. London Math. Soc. 47 (1993)\textcolor[rgb]{1,0,0}{,} 227-239.
\bibitem{P4} C.E. Praeger, Finite transitive permutation groups and finite vertex-transitive graphs, in: Graph Symmetry: Algebraic Methods and Applications, in: NATO Adv. Sci. Inst. Ser. C, 497 (1997)\textcolor[rgb]{1,0,0}{,} 277-318.
\bibitem{P5} C.E. Praeger, Imprimitive symmetric graphs, Ars Combin. A 19 (1985)\textcolor[rgb]{1,0,0}{,} 149-163.
\bibitem{Sabidussi}
B. O. Sabidussi, Vertex-transitive graphs, Monash Math. 68 (1964), 426-438.
\bibitem{Schur}
J. Schur, \" Uber die Darstellung der endlichen Gruppen durch gebrochene
lineare Substitutionen, J. Reine Angew. Math. 127 (1904), 20-50.
\bibitem{Wang}
Y. Wang, Y.-Q. Feng, J.-X. Zhou, Automorphism of Cayley digraph of $2$-genetic groups of prime-power order,
J. Combin Theory Ser. A, 143 (2016), 88-106.
\bibitem{Weiss}
R.M. Weiss, Presentations for $(G,s)$-transitive graphs of small valency,
Math. Proc. Cambridge Philos. Soc. 101 (1987), 7-20.
\bibitem{Wilson}
R. A. Wilson, The Finite Simple Groups, Springer-Verlag, London, 2009.
\bibitem{XFWX2005}
S. J. Xu, X. G. Fang, J. Wang, and M. Y. Xu, On cubic $s$-arc transitive Cayley graphs of finite simple groups,
Europ. J. Combin. 26 (2005), 133-143.
\bibitem{XFWX}
S. J. Xu, X. G. Fang, J. Wang, and M. Y. Xu, $5$-arc transitive cubic Cayley graphs on finite simple groups,
Europ. J. Combin. 28 (2007), 1023-1036.
\end{thebibliography}
\end{document}